# Two Simple Ways of Generating the Partitions of (n+1) from the Partitions of n.

Dhananjay P. Mehendale
Department of Electronics, S. P. College, Pune-411009, India.

**1. Introduction:**

A partition of a nonnegative integer n is a representation of n as a sum of positive integers called parts of the partition. The order of parts is irrelevant, and so partitions are generally written in non-increasing order, thus, the partitions of 5 are:
5, 4+1, 3+2, 3+1+1, 2+2+1, 2+1+1+1, 1+1+1+1+1.
We denote by P(n) the number of partitions of n, thus P(5) = 7.
The famous generating function for partitions is due to Euler [1]. It is stated as:

$$\sum_{n=o}^{\infty} P(n)q^n = \prod_{j=1}^{\infty} \frac{1}{1-q^j}, where\ |q|<1 \cdots (1).$$

So, by expanding the infinite products as power series, carrying out multiplication of these series, and then comparing the coefficients of $q^n$ one can generate as well as count the partitions of n.
In this paper I propose two new methods to generate the partitions of (n+1) from the partitions of n.
The first method that I have proposed in this short paper is based on the observation that among the partitions of n there are two kinds of partitions, say partitions of the "First Kind", having last two (or more) parts identical, and partitions of the "Second Kind", having only one part or having last two parts strictly different. We see that partitions of the "First Kind" generate one and only "one" partition of (n+1) while partitions of the "Second Kind" generate "two" partitions of (n+1).
As n is made large the count of partitions of the "Second Kind" becomes large. This fact leads to the astronomical growth of the number of partitions of n, P(n), with the increase in n.

**2. The First Method:**

This method is based on the following simple observations:

Among the partitions of n,
(a) Partitions of n of the "First Kind" produce exactly one partition of (n+1), obtained by simply adding "1" separately.
(b) Partitions of n of the "Second Kind" produce exactly two partitions of (n+1), one obtained as above by adding "1" separately, while the other obtained by augmenting 1 in the last part (or in the only part), i.e. if the last part is $\alpha_k$ then one has to make it "$(\alpha_k+1)$" by augmentation.

**Lemma 1:** For *every* partition of (n+1) the method implies a *unique* predecessor in the partitions of n.

**Proof:** Let $\alpha$ be a k- partition (partition having k parts) of (n+1),

$$\alpha = \alpha_1 + \alpha_2 + \alpha_3 + \cdots + \alpha_{(k-1)} + \alpha_k$$

(i) If $\alpha_k = 1$ then $\alpha$ is uniquely obtained, as per the method, from the (k-1)-partition of n, $\beta$ say, by adding "1" separately where

$$\beta = \alpha_1 + \alpha_2 + \alpha_3 + \cdots + \alpha_{(k-1)}$$

(ii) If $\alpha_k > 1$ then again $\alpha$ is uniquely obtained from the k-partition of n, $\gamma$ say, by augmenting "1" in the last part where

$$\gamma = \alpha_1 + \alpha_2 + \alpha_3 + \cdots + \alpha_{(k-1)} + (\alpha_k - 1).$$

Thus, every partition of (n+1) is uniquely derivable from a partition of n. The word "every" ensures that all the partitions (without omission) of (n+1) are achievable by this method, and the word "uniquely" implies that every partition of (n+1) gets generated only once (without repetition). Hence the lemma.

Algorithm 1:

(1) Make two groups of the two kinds of the partitions of n:
   (i) A group of the partitions of n with $\alpha_k = \alpha_{(k-1)}$, and
   (ii) A group with k=1 or, when k > 1, with $\alpha_k < \alpha_{(k-1)}$, for all k > 1.

(2) From each partition of n in the first group obtain a partition of (n+1) by separately adding a unit i.e. "1" in each partition.

(3) From each partition of n in the second group obtain exactly two partitions of (n+1), one by adding "1" separately (as in (2)) and the other by augmenting "1" in the last part i.e. by changing $\alpha_k$ to $(\alpha_k+1)$.

**Theorem:** Algorithm 1 generates all the partitions of (n+1) from the partitions of n, without omission or repetition.

**Proof:** Follows from the lemma 1.

**Example:** Let n = 5. Referring to the list of the partitions of 5 given above, the groups, obtained as per step (1) of the algorithm 1, are:
Group 1: {3+1+1, 2+1+1+1, 1+1+1+1+1}.
Group 2: {5, 4+1, 3+2, 2+2+1}
From the first group, the partitions of 6, obtained as per step (2) of the algorithm 1, are:
  {3+1+1+1, 2+1+1+1+1, 1+1+1+1+1+1}.
From second group, the partitions of 6 that we get, as per step (3) of the algorithm 1, first by adding "1" separately and secondly by augmentation of "1" in the last or the only part, are:
{5+1, 4+1+1, 3+2+1, 2+2+1+1, 6, 4+2, 3+3, 2+2+2}.
We are thus ready with the complete list of the partitions of 6.

Let us denote by Q(n) the number of partitions of n containing k parts, with k = 1, and the partitions of n with last two parts strictly different, i.e. $\alpha_k < \alpha_{(k-1)}$, when k > 1.

Then we have, as an immediate consequence of the above algorithm, the following

**Theorem:** P(n+1) = P(n) + Q(n).

**Proof:** (i) P(n) on the r.h.s. counts the partitions of (n+1) obtained by adding "1" separately, in all the partitions of n.
  (ii) Q(n) on the r.h.s. counts the partitions of (n+1) obtained from the partitions of n having the only part or having the last two parts strictly different, by augmenting "1" in the only part i.e. n, or, by augmenting "1" in the last part where the last two parts strictly differ.

**A generating function for Q(n):**

**Theorem:** The following is the generating function for Q(n), the partitions of n having the only part i.e. n, or having strictly different last two parts:

$$\sum_{n=o}^{\infty} Q(n)q^n = \sum_{s=1}^{\infty} q^s \left[ \prod_{j=(s+1)}^{\infty} \left( \frac{1}{1-q^j} \right) \right] \cdots (2)$$

**Proof:** Expanding each part in the product on the right side of the Euler's generating function given in Equation.(1) above we have

$$\sum_{n=0}^{\infty} P(n)q^n = \prod_{j=1}^{\infty} \left( 1 + q^j + q^{2 \cdot j} + q^{3 \cdot j} + \cdots + q^{m \cdot j} + \cdots \right)$$

after rearranging the right side the above equation can be rewritten as

$$\sum_{n=0}^{\infty} P(n)q^n = \sum_{k=1}^{\infty} q^k \prod_{j=(k+1)}^{\infty} \left( 1 + q^j + q^{2 \cdot j} + q^{3 \cdot j} + \cdots + q^{m \cdot j} + \cdots \right)$$

+ other terms.

It is easy to check that in the sum of products on the r.h.s. of the above equation a product for (any) k counts the partitions with the *only* or the *last* part equal to k while all other parts *strictly bigger* than k. So, by collecting all terms whose indices add to n, for all k, we can obtain all the partitions of n that are of the "Second Kind".

It is straightforward to see that, after rewriting the sum of the products in the above equation, the equation can be rewritten as

$$\sum_{n=0}^{\infty} P(n)q^n = \sum_{k=1}^{\infty} q^k \prod_{j=(k+1)}^{\infty} \left( \frac{1}{1-q^j} \right)$$

+ other terms.

Hence the result.

### 3. The Second Method:

This method is similar to the first method. The difference is in the definition for classification of the partitions of n.

Consider a partition of n:

$$\alpha = \alpha_1 + \alpha_2 + \alpha_3 + \cdots + \alpha_k + 1 + 1 + \cdots + 1,$$

where $\aleph$ are the number of units present in the partition.
Hereafter, a partition of n will be called the partition of the "First Kind" if $\aleph = 0$, or $\aleph \geq \alpha_k$, and a partition of n will be called the partition of the "Second Kind" when $1 \leq \aleph < \alpha_k$.

In this method,
 (a) Partitions of n of the "First Kind" produce exactly one partition of (n+1), obtained by simply adding "1" separately.
 (b) Partitions of n of the "Second Kind" produce exactly two partitions of (n+1), one obtained as above by adding "1" separately, while the other obtained by augmenting 1 in the count of the number of units i.e. replacing all units, $\aleph$ in number, by a single number ($\aleph+1$). Note that in this procedure of generating the partitions of (n +1) from the partitions of n we are not augmenting "1" to the last part but collecting the units present in the partition (of n) +1 into a single number. So, the partitions of n with $\aleph = 0$ produce only one partition, produced by adding "1" separately. Therefore, in this method the partition of (n+1) with only one part i.e. {n+1} is not generated and has to be added in the list separately.

**Lemma 2:** For *every* partition of (n+1) (except {n+1}) the method (second) implies a *unique* predecessor in the partitions of n.

**Proof:** Let $\alpha$ be a k- partition (partition having k parts) of (n+1), and suppose k > 1:

$$\alpha = \alpha_1 + \alpha_2 + \alpha_3 + \cdots + \alpha_{(k-1)} + \alpha_k$$

(i) If $\alpha_k = 1$ then $\alpha$ is uniquely obtained, as per the method, from the
(k-1)-partition of n, $\beta$ say, by adding "1" separately where

$$\beta = \alpha_1 + \alpha_2 + \alpha_3 + \cdots + \alpha_{(k-1)}$$

(ii) If $\alpha_k > 1$ then again $\alpha$ is uniquely obtained from the $(k+\alpha_k-2)$-partition of n, $\gamma$ say, by collecting $(\alpha_k-1)$ units into a single number and augmenting an additional "1" in it to make it $\alpha_k$, where

$$\gamma = \alpha_1 + \alpha_2 + \alpha_3 + \cdots + \alpha_{(k-1)} + 1 + 1 + \cdots + 1,$$

with $(\alpha_k - 1)$ units.

Thus, every partition of (n+1) (except the partition with one part, i.e. {n+1}) is uniquely derivable from a partition of n. The word "every" ensures that all the partitions (without omission) of (n+1) are achievable by this method, and the word "uniquely" implies that every partition of (n+1) gets generated only once (without repetition).
Hence the lemma.

Algorithm 2:

From the given partitions of n having k parts,
(1) Make two groups of the two kinds of the partitions of n:
(i) A group with $\aleph = 0$, or with $\aleph \geq \alpha_{(k-1)}$, and
(ii) A group with $1 \leq \aleph < \alpha_{(k-1)}$, for all k > 1,

where $\aleph$ denotes the number of units present in the partition..
(2) From each partition of n in the first group obtain a partition of (n+1) by separately adding a unit i.e. "1" in each partition.
(3) From each partition of n in the second group obtain exactly two partitions of (n+1), one by adding "1" separately (as in (2)) and the other by replacing the separate $\aleph$ number of units by a single number $(\aleph+1)$
(4) Add the partition {n+1} separately in the list.

**Theorem:** Algorithm 2 generates all the partitions of (n+1) from the partitions of n, without omission or repetition.

**Proof:** Follows from lemma 2, and by the action of separately adding the partition of (n+1) containing only one part i.e. {n+1}.

**Example:** Let n = 5. Referring to the list of the partitions of 5 given above, the groups, as per the step (1) of the algorithm 2, are:
Group 1: {5, 3+2, 2+1+1+1, 1+1+1+1+1}.
Group 2: {4+1, 3+1+1, 2+2+1}
From the first group, the partitions of 6 obtained as per step (2) of the algorithm 2, are:

{5+1, 3+2+1, 2+1+1+1+1, 1+1+1+1+1+1}.

From second group, the partitions of 6 that we get as per step (3) of the algorithm 2, first by adding "1" separately and secondly by augmentation of "1" in the last part, are:

{4+1+1, 3+1+1+1, 2+2+1+1, 4+2, 3+3, 2+2+2}.

Adding to the list, as per step (4) of the algorithm 2, we get the partition {6}

We are thus ready with the complete list of the partitions of 6.

**4. Conclusion:** Looking closely at the partitions of the successive numbers one has at his /her service two simple and interesting recursive procedures for the generation of partitions.

## Acknowledgements

The author is thankful to Dr. M. R. Modak and Dr. S. A. Katre, Bhaskaracharya Pratishthana, Pune, for their keen interest